\documentclass[12pt]{article}
\usepackage{hyperref}
\usepackage{amsfonts}
\newtheorem{theorem}{Theorem}

\newtheorem{example}[theorem]{Example}
\newtheorem{definition}[theorem]{Definition}
\newcommand{\R}{\mathbb{R}}
\newcommand{\grad}{\mbox{\rm grad}}
\newcommand{\gradl}{\mbox{\rm grad}^L}
\newcommand{\gradn}{\mbox{\rm grad}^N}
\newcommand{\gradm}{\mbox{\rm grad}^M}
\newcommand{\gradp}{\mbox{\rm grad}^P}
\newcommand{\hess}{\mbox{Hess }}

\newcommand{\po}{{\hspace*{-1ex}}{\bf .  }}
\newcommand{\ii}{isometric immersion }
\def\<{\langle}

\def\nl{\nabla^L}
\def\np{\nabla^P}
\def\>{\rangle}
\def\bea{\begin{eqnarray*} }
\def\eea{\end{eqnarray*} }
\def\be{\begin{equation} }
\def\ee{\end{equation} }
\def\nab{\widetilde\nabla}

\begin{document}

\title{A mean curvature estimate for cylindrically bounded submanifolds}
\author {Luis J. Al\'{i}as and Marcos Dajczer}
\date{}
\maketitle

\begin{abstract} We extend the estimate obtained in \cite{AlBeDa}
for the mean curvature of a cylindrically bounded proper submanifold
in a product manifold with  an Euclidean space as one factor
to a general product ambient space endowed with a warped product structure.
\end{abstract}

Let ($L^\ell,g_L)$ and $(P^n,g_P)$ be complete Riemannian manifolds of
dimension $\ell$ and $n$, respectively, where $L^\ell$ is non compact. Then, let
$N^{n+\ell}=L^\ell\times_{\rho} P^n$
be the product manifold $L^\ell\times P^n$
endowed with the warped product metric $ds^2=dg_L+\rho^2dg_P$
for some positive warping function $\rho\in C^\infty(L)$.

Let $B_P(r_0)$ denote the geodesic ball  with radius $r_0$ centered  at a
reference point $o\in P^n$. We assume that the radial sectional curvatures
in $B_P(r_0)$ along the geodesics issuing from $o$  are bounded as
$K_P^{\mathrm{rad}}\leq b$ for some constant $b\in\R$,
and that $0<r_0<\min\{{\rm inj}_P(o), \pi/2\sqrt{b}\}$ where ${\rm inj}_P(o)$
is the injectivity radius at $o$ and $\pi/2\sqrt{b}$ is replaced  by $+\infty$
if $b\leq 0$. Then, the mean curvature of the geodesic sphere
$S_P(r_0)=\partial B_P(r_0)$ can be estimated from below by the mean curvature of
a geodesic sphere of a space form of curvature $b$, namely,
$$
C_b(t)
=\left\{\begin{array}{lll}
\sqrt{b}\cot(\sqrt{b}\, t) & \mathrm{if} & b >0,\\
1/t & \mathrm{if} & b =0,\\
\sqrt{-b}\coth(\sqrt{-b}\, t) & \mathrm{if}  & b <0.
\end{array}\right.
$$
This is a direct consequence of the comparison theorems for the Riemannian distance, since the Hessian (respectively, Laplacian)
of the distance function is nothing but the second fundamental form (respectively, mean curvature) of the geodesic spheres.
A classical reference about this topic is \cite{GW}. We also refer the reader to \cite{Petersen} or \cite{PRSbook} for
a modern approach to the Hessian and Laplacian comparison theorems.

By a \emph{cylinder} in the warped space $N^{n+\ell}$ we mean a closed subset of the form
$$
{\cal C}_{r_0}=\{(x,y)\in N^{n+\ell}:x\in L^\ell\; \mbox{and} \;y\in B_P(r_0)\}.
$$
Since the submanifolds $L^\ell\times\{p_0\}\subset N^{n+\ell}$ are totally geodesic,
we have that
$$
|\rho H_{{\cal C}_{r_0}}|\geq \frac{n-1}{\ell+n-1}C_b(r_0)
$$
where $H_{{\cal C}_{r_0}}$ is the mean curvature vector field  of the hypersurface
$L^\ell\times S_p(r_0)$.

The following theorem extends the result in \cite{AlBeDa} where the cylinders under
consideration are contained in product spaces  $\R^\ell\times P^n$.
After the statement,  we recall  from \cite{AlBeMoPi}  the concept of an Omori-Yau pair
on a Riemannian manifold and discuss some implications of its existence.

\begin{theorem}\label{main}\po
Let $f\colon\,M^m\to L^\ell\times_\rho P^n$ be an \ii where~$L^\ell$ carries an
Omori-Yau pair for the Hessian and the functions $\rho$ and $|\grad\log\rho|$ are bounded.
If $f$ is proper and  $f(M)\subset {\cal C}_{r_0}$, then
$\sup_M|H|=+\infty$ or
\be\label{four}
\sup_M\rho|H|\geq\frac{m-\ell}{m}C_b(r_0)
\ee
where $H$ is the mean curvature vector field of $f$.
\end{theorem}

In the proof we see that the existence in $L^\ell$ of a Omori-Yau pair for the Hessian
provides conditions, in a  function theoretic form,
that guarantee the validity of the Omori-Yau Maximum Principle on $M^m$
in terms of the corresponding property of $L^\ell$ and the geometry of the immersion.

\begin{definition}\po\label{def} {\em The pair of functions $(h,\gamma)$ for $h\colon\,\R_+\to\R_+$ and
$\gamma\colon\, M\to\R_+$ form an \emph{Omori-Yau pair for the Hessian}
in $M$ if they  satisfy:
\begin{enumerate}
\item[(a)]
$h(0)>0$ and $h'(t)\geq 0$ for all $t\in \R_+$,
\item[(b)]
$\limsup\limits_{t\to+\infty}\,\displaystyle t h\big(\sqrt t\big)/h(t)<+\infty$,
\item[(c)]
$\displaystyle\int_0^{+\infty}\mathrm dt/\sqrt{\strut h(t)}=+\infty$,
\item[(d)] The function $\gamma$ is proper,
\item[(e)] $|\grad\,\gamma|\le c\sqrt\gamma$\, for some $c>0$
outside a compact subset of $M$,
\item[(f)]  $\hess\gamma\le d\sqrt{\strut\gamma h(\sqrt\gamma)}$\, for some
$d>0$ outside a compact subset of $M$.
\end{enumerate}
}\end{definition}
Similarly, the  pair $(h,\gamma)$ forms an \emph{Omori-Yau pair for the Laplacian} in $M$
if they satisfy conditions $(a)$ to $(e)$  and
\begin{enumerate}
\item[(f')] $\Delta\gamma\le d\sqrt{\strut\gamma h(\sqrt\gamma)}$ for some
$d>0$ outside a compact subset of $M$.
\end{enumerate}
\medskip

The following fundamental result due to Pigola, Rigoli and Setti \cite{PiRiSe} gives
sufficient conditions  for an Omori-Yau Maximum Principle
to hold for a Riemannian manifold.

\begin{theorem}\po\label{sufcondOY}   Assume that a Riemannian manifold $M$
carries  an \mbox{Omori-Yau} pair for the Hessian (respec., Laplacian).
Then, the Omori--Yau Maximum Principle for the Hessian (respec., Laplacian) holds in $M$.
\end{theorem}

Recall that the \emph{Omori--Yau Maximum Principle for the Hessian} holds
for $M$ if for any function $g\in C^\infty(M)$ bounded
from above there exists a sequence of points $\{p_k\}_{k\in\mathbb N}$ in $M$ such that
\begin{itemize}
\item[(a)] $\lim\limits_{k\to\infty}g(p_k)=\sup\limits_Mg$,
\item[(b)] $|\grad\, g(p_k)|\le 1/k$,
\item[(c)]  $\hess g(p_k)(X,X)\le (1/k)g_M(X,X)$ for
all $X\in T_{p_k}M$.
\end{itemize}
Similarly, the \emph{Omori--Yau Maximum Principle for the Laplacian} holds
for $M$ if the above properties are satisfied with $(c)$ replaced by
\begin{itemize}
\item[(c')] $\Delta g(p_k)\le 1/k$.
\end{itemize}

\begin{example}\po\label{two}
{\em Let $M^m$ be a complete but non compact Riemannian manifold and denote $r(y)=\mbox{dist}_M(y,o)$
for some reference point $o\in M^m$.  Assume that the radial sectional curvature of $M^m$ satisfies
$K^{\mathrm{rad}}\geq -h(r)$, where the smooth function $h$ satisfies $(a)$ to $(c)$ in
Definition \ref{def} and is  even at the origin, that is, $h^{(2k+1)}(0)=0$ for $k\in {\mathbb N}$.
Then, it was shown in \cite{PiRiSe}  that the functions $(h,r^2)$ form an Omori-Yau pair for the
Hessian. As for the function $h$, one can choose
$$
h(t)=t^2\prod_{j=1}^N(\log^{(j)}(t))^2,\;\;t\gg 1,
$$
where $\log^{(j)}$ stands for the $j$-th iterated logarithm.
}\end{example}

To conclude this section, we first observe that Theorem \ref{main} is sharp.  This is clear from
(\ref{four}) by taking as $P^n$  a space-form and as $M$ the hypersurface $L^\ell\times S_P(r_0)$
in $N^{n+\ell}$. Moreover, in view of Example \ref{two} it follows taking $L^\ell=\R^\ell$ and
constant $\rho$ that we recover the result in \cite{AlBeDa}.

\section{The proof}

We first introduce some additional notations and then we recall a few basic facts on warped product
manifolds.
\medskip

Let $\<\,,\,\>$ denote the metrics in $N^{n+\ell}$, $L^\ell$ and $M^m$ whereas
$(\,,\,)$ stands for the metric in $P^n$. The corresponding norms  are  $|\;\;|$ and $\|\;\;\|$.
In addition, let $\nabla$ and $\nab$  denote the Levi-Civita connections in $M^m$ and $N^{n+\ell}$,
respectively, and  $\nl$  and $\np$ the ones in $L^\ell$ and $P^n$.

We always denote vector fields in $TL$ by $T,S$ and in $TP$ by $X,Y$. In addition, we identify vector fields in $TL$
and $TP$ with \emph{basic} vector fields in $TN$  by taking $T(x,y)=T(x)$ and $X(x,y)=X(y)$.

For the Lie-brackets of basic vector fields, we have that $[T,S]\in TL$ and $[X,Y]\in TP$ are basic
and that $[X,T]=0$.  Then, we have
$$
\nab_ST=\nl_ST,
$$
$$
\nab_XT=\nab_TX=T(\varrho)X
$$
and
$$
\nab_XY=\np_XY-\<X,Y\>\gradl\varrho
$$
where the vector fields $X,Y$ and $T$ are basic and  $\varrho=\log\rho$.
\medskip

Our proof follows the main steps in \cite{AlBeMoPi}, where the geometric situation considered differs from ours since
there $f(M)$ is contained in a \textit{cylinder} of the form
\[
\{ (x,y)\in N^{n+\ell} : x\in B_L(r_0)\; \mbox{and} \; y\in P^n\}.
\]
In fact, a substantial part of the argument is to show
that the Omori-Yau pair for the Hessian in $L^\ell$ induces an Omori-Yau pair
for the Laplacian for a non compact $M^m$ when $|H|$ is bounded. Thus, the  Omori--Yau Maximum Principle
for the  Laplacian holds in $M^m$, and the proof follows from a application of the latter.
\medskip

Suppose that $M^m$ is non compact and let $(h,\Gamma)$ be an Omori-Yau  pair for the Hessian in $L^\ell$.
For $p\in M^m$ denote $f(p)=(x(p),y(p))$.
Set  $\tilde\Gamma(x,y)=\Gamma(x)$ for $(x,y)\in N^{n+\ell}$ and
$$
\gamma(p)=\tilde\Gamma(f(p))=\Gamma(x(p)).
$$

We show next that $(h,\gamma)$ is an Omori-Yau pair for the Laplacian in $M^m$.
First, we argue that the  function $\gamma$ is proper. To see this, let $p_k\in M^m$ be a divergent sequence, i.e.,
$p_k\to\infty$ in $M^m$ as $k\to +\infty$. Thus, $f(p_k)\to \infty$ in $N^{n+\ell}$ since $f$ is proper.
Since $f(M)$ lies inside a cylinder, then $x(p_k)\to \infty$
in $L^\ell$. Hence, $\gamma(p_k)\to+\infty$ as $k\to +\infty$  since $\Gamma$ is proper, and thus $\gamma$ is proper.

It remains to verify conditions  $(e)$ and
$(f')$ in Definition \ref{def}.
We have from $\tilde\Gamma(x,y)=\Gamma(x)$ that
$$
\<\gradn\tilde\Gamma(x,y),X\>=0.
$$
Thus,
$$
\gradn\tilde\Gamma(x,y)=\gradl \Gamma(x).
$$
Since $\gamma=\tilde\Gamma\circ f$, we obtain
\be\label{six}
\gradn\tilde \Gamma(f(p))=\gradm \gamma(p) + \gradn\tilde\Gamma(f(p))^\perp
\ee
where $(\;)^\perp$ denotes taking the normal component to $f$. Then,
$$
|\gradm\gamma(p)|\leq |\gradn \tilde\Gamma(f(p))|=|\gradl\Gamma(x(p))|
\leq c\sqrt{\Gamma(x(p))}=c\sqrt{\gamma(p)}
$$
outside a compact subset of $M^m$, and thus $(e)$ holds.

We have that
$$
\nab_T\gradn\tilde\Gamma=\nl_T\gradl\Gamma.
$$
Hence,
$$
\hess\tilde\Gamma(T,S)=\hess\Gamma(T,S)
$$
and
$$
\hess\tilde\Gamma(T,X)=0.
$$
Moreover,
$$
\nab_{X}\gradn\tilde\Gamma=\nab_{X}\gradl\Gamma
=\gradl\Gamma(\varrho)X.
$$
Hence,
$$
\hess\tilde\Gamma(X,Y)=\<\gradl\Gamma,\gradl\varrho\>\<X,Y\>.
$$
For a unit vector $e\in T_pM$, set $e=e^L+e^P$ where $e^L\in T_{x(p)}L$ and $e^P\in T_{y(p)}P$.  Then,
$$
\hess\tilde\Gamma(f(p))(e,e)
=\hess\Gamma(x(p))(e^L,e^L)+ \<\gradl\Gamma(x(p)),\gradl\varrho(x(p))\>|e^P|^2.
$$
Moreover, an easy computation using (\ref{six}) yields
$$
\hess\gamma(p)(e,e)=\hess\tilde\Gamma(f(p))(e,e) + \<\gradl\Gamma(x(p)),\alpha(p)(e,e)\>
$$
where $\alpha$ denotes the  second fundamental of $f$ with values in the normal bundle.
Thus,
\bea
\hess\gamma(p)(e,e) \!\!\!&=& \!\!\! \hess\Gamma(x(p))(e^L,e^L)+ \<\gradl\Gamma(x(p)),\gradl\varrho(x(p))\>|e^P|^2 \\
 \!\!\!&+&\!\!\!\<\gradl\Gamma(x(p)),\alpha(p)(e,e)\>.
\eea
Since $\hess\Gamma\leq d\sqrt{\Gamma h(\sqrt{\Gamma})}$ for some positive constant $d$ outside a compact subset of
$L^\ell$ and the immersion is proper, then
$$
\hess\Gamma(x(p))(e^L,e^L)\leq d\sqrt{\gamma(p)h(\sqrt{\gamma(p)})}|e^L|^2\leq d\sqrt{\gamma(p)h(\sqrt{\gamma(p)})}
$$
outside a compact subset of $M^m$. From
$|\gradl\Gamma|\leq c\sqrt{\Gamma h(\sqrt{\Gamma})}$ for some $c$ outside a compact subset of
$L^\ell$ and $\sup_L|\gradl\varrho|<+\infty$, we  have
$$
\<\gradl\Gamma(x(p)),\gradl\varrho(x(p))\>|e^P|^2\leq c'\sqrt{\gamma(p)}
$$
for some positive constant $c'$ outside a compact subset of $M^m$. Being $\gamma$ proper and $h$
unbounded from $(a)$ and $(b)$ in Definition
\ref{def}, then
$$
\sqrt{\gamma}\leq \sqrt{\gamma h(\sqrt{\gamma})}
$$
outside a compact subset of $M^m$, since $\gamma\rightarrow+\infty$ as $p\rightarrow\infty$ and
$\lim_{t\rightarrow+\infty}h(t)=+\infty$.
Thus, we obtain
\be\label{nine}
\hess\gamma(e,e)\leq d_1\sqrt{\gamma h(\sqrt{\gamma})}+\<\gradl\Gamma(x),\alpha(e,e)\>
\ee
for same constant $d_1>0$, outside a compact subset of $M^m$.

On the other hand, we may assume that
\be\label{three}
|H|\le c\,\sqrt{\strut h(\sqrt\gamma)}
\ee
for some constant $c>0$, outside a compact subset of $M^m$.
Otherwise, there exists a sequence
$\{p_k\}_{k\in\mathbb{N}}$ in $M^m$ such that $p_k\to \infty$ as $k\to +\infty$ and
$$
|H(p_k)|>k\sqrt{\strut h(\sqrt{\gamma(p_k)})}.
$$
Being $\gamma$ proper and $h$ unbounded from $(a)$ and $(b)$ in Definition \ref{def},
we conclude that $\sup_M|H|=+\infty$, in which case we are done with the proof of the theorem.

We obtain from (\ref{nine}) using (\ref{three}) that
$$
\Delta\gamma\leq c_1\sqrt{\gamma h(\sqrt{\gamma})}
$$
for some constant $c_1>0$ outside a compact subset of $M^m$, and thus $(f')$ has been proved.

Consider the distance function $r(y)=\mbox{dist}_P(y,o)$ in $B_P(r_0)$ and
define $\tilde r\in  C^\infty(N)$ by $\tilde r(x,y)=r(y)$. Then,
$$
\<\gradn\tilde r(x,y),T\>=0.
$$
Thus,
$$
\rho^2(x)\gradn\tilde r(x,y)=\gradp r(y).
$$
We obtain that
$$
\nab_T\gradn\tilde r=\nab_T(\rho^{-2}\gradp r)=-\rho^{-2}T(\varrho)\gradp r.
$$
Therefore,
$$
\hess\tilde r(T,S)=0
$$
and
$$
\hess\tilde r(T,X)=-\rho^{-2}T(\varrho)\<\gradp r,X\>= -T(\varrho)(\gradp r,X).
$$
Moreover,
$$
\nab_{X}\gradn\tilde r=\nab_{X}(\rho^{-2}\gradp r)
=\rho^{-2}\left(\np_X\gradp r-\<X,\gradp r\>\gradl\varrho\right).
$$
Hence,
$$
\hess\tilde r(X,Y)=\rho^{-2}\<\np_{X}\gradp r,Y\>
=(\np_X\gradp r,Y)=\hess r(X,Y).
$$
For $e\in TM$, we have
$$
\hess\tilde r(e,e)=-2\<\gradl\varrho,e\>(\gradp r,e^P)+\hess r(e^P,e^P).
$$
From the Hessian comparison theorem (cf. Chapter 2 in \cite{PRSbook} for a modern approach), we obtain
$$
\hess r(e^P,e^P)\geq C_b(r)(\|e^P\|^2 - (\gradp r,e^P)^2).
$$
Therefore,
\be\label{ten}
\hess\tilde r(e,e)\geq -2\<\gradl\varrho,e\>(\gradp r,e^P)+C_b(r)(\|e^P\|^2 - (\gradp r,e^P)^2).
\ee

We define
$u\in C^\infty(M)$ by
$$
u(p)=r(y(p)).
$$
Thus,
$u=\tilde r\circ f$ and
\be\label{one}
\gradn\tilde r(f(p))=\gradm u(p) + \gradn\tilde r(f(p))^\perp.
\ee
Using (\ref{one}) gives
$$
\hess u (e_i,e_j) = \hess \tilde r(e_i,e_j) + \<\gradn\tilde r, \alpha(e_i,e_j)\>
$$
where  $e_1,\ldots,e_m$ an orthonormal frame of $TM$. Thus,
\be\label{twelve}
\Delta u= \sum_{j=1}^m\hess\tilde r(e_j,e_j) + m\<\gradn\tilde r, H\>.
\ee

We have from $e_j=e_j^L+e^P_j$ that
$$
1=\<e_j,e_j\>=\rho^2\|e_j^P\|^2 + \sum_{k=1}^\ell\<e_j,T_k\>^2
$$
where $T_1,\ldots,T_\ell$ is an orthonormal frame for $TL$. Hence,
$$
m=\rho^2\sum_{j=1}^m\|e_j^P\|^2 +\sum_{k=1}^\ell|T_k^\top|^2,
$$
where $T^\top$ is the tangent component of $T$.  We obtain that
\be\label{seven}
\sum_{j=1}^m\|e_j^P\|^2\geq (m-\ell)\rho^{-2}.
\ee

We obtain from (\ref{ten}) and
$$
(\gradp r,e_j^P)=\<\gradn\tilde r,e_j^P\>=\<\gradn\tilde r,e_j\>=\<\gradm u,e_j\>
$$
that
$$
\hess\tilde r(e_j,e_j) \geq-2\<\gradl\varrho,e_j\>\<\gradm u,e_j\>
+ C_b(u)(\|e_j^P\|^2 - \<\gradm u,e_j\>^2).
$$
Taking trace and using (\ref{seven}) gives
$$
\sum_{j=1}^m\hess\tilde r(e_j,e_j) \geq-2\<\gradl\varrho,\gradm u\>
+ C_b(u)\left((m-\ell)\rho^{-2}-|\gradm u|^2\right).
$$
Since
$$
\<\gradn\tilde r,\gradn\tilde r\>=\rho^2(\rho^{-2}\gradp r,\rho^{-2}\gradp r)=\rho^{-2},
$$
we have
$$
\<\gradn\tilde r, H\>\geq - \rho^{-1}|H|.
$$
We conclude using (\ref{twelve}) that
$$
\Delta u\geq-2\<\gradl\varrho,\gradm u\>
+ C_b(u)\left((m-\ell)\rho^{-2}-|\gradm u|^2\right)-m\rho^{-1}|H|.
$$
Thus,
$$
\rho|H|\geq \frac{m-\ell}{m}C_b(u)-\frac{\rho^2}{m}\left(\Delta u+2|\gradl\varrho||\gradm u|
+C_b(u)|\gradm u|^2 \right).
$$

If $M^m$ is compact, the proof follows easily by computing the inequality  at a point
of maximum of $u$. Thus, we may now assume that $M^m$ is non compact and that (\ref{three})
holds.

Since $f(M)\subset {\cal C}_{r_0}$, we have  $u^*=\sup_Mu\leq r_0<+\infty$.
By the Omori-Yau maximum principle there is a sequence $\{p_k\}_{k\in\mathbb{N}}$ in  $M^m$
such that
$$
u(p_k)>u^*-1/k,\;\;\;|\gradm u(p_k)|<1/k\;\;\mbox{and}\;\;\Delta u(p_k)<1/k.
$$
By assumption, we have $\sup_L\rho=K_1<+\infty$ and $\sup_L|\gradl\varrho|=K_2<+\infty$. Hence,
$$
\sup_M\rho|H|\geq \rho(p_k)|H(p_k)|\geq
\frac{m-\ell}{m}C_b(u(p_k))-\frac{K_1^2}{m}\left(\frac{1+2K_2}{k}+\frac{1}{k^2}C_b(u(p_k))\right).
$$
Letting $k\to +\infty$, we obtain
$$
\sup_M\rho|H|\geq \frac{m-\ell}{m}C_b(u^*)\geq \frac{m-\ell}{m}C_b(r_0),
$$
and this concludes the proof of the theorem.

\vspace*{-1ex}

{\renewcommand{\baselinestretch}{1}
\hspace*{-20ex}\begin{tabbing}
\indent \=Luis J. Alias \hspace{25,5ex} Marcos Dajczer  \\
\>Departamento de Matematicas  \hspace{7,5ex} IMPA\\
\>Universidad de Murcia \hspace{15,4ex}
Estrada Dona Castorina, 110\\
\>  Campus de Espinardo
\hspace{17ex}22460-320 --- Rio de Janeiro ---RJ\\
\> E-30100 Espinardo, Murcia\hspace{12,6ex}Brazil \\
\> Spain\hspace{33,6ex} marcos@impa.br\\
\> ljalias@um.es
\end{tabbing}}


\begin{thebibliography}{20}

\bibitem{AlBeDa} L. J. Al\'{i}as, G. P. Bessa, M. Dajczer, \emph{The mean curvature of
cylindrically  bounded submanifolds}, Math.\ Ann.\ \textbf{345}, (2009) 367--376.

\bibitem{AlBeMoPi} L. J. Al\'\i as, G. P. Bessa, J. F. Montenegro and P. Piccione.
\emph{Curvature estimates for submanifolds in warped products.} To appear in Results in Mathematics.
Available at http://arxiv.org/abs/1009.3467

\bibitem{GW} R.E. Greene, H. Wu, Function theory on manifolds which possess a pole. Lecture Notes in Mathematics, 699.
Springer, Berlin, 1979.

\bibitem{Petersen} P. Petersen, Riemannian geometry. Second edition. Graduate Texts in Mathematics, 171. Springer,
New York, 2006.

\bibitem{PiRiSe} S. Pigola, M. Rigoli, A. Setti, \emph{Maximum Principle on Riemannian
Manifolds and Applications}, Memoirs Amer.\ Math.\ Soc.\ \textbf{822} (2005).

\bibitem{PRSbook} S. Pigola, M. Rigoli, A. Setti, {Vanishing and finiteness results in geometric analysis. A
generalization of the Bochner technique}. Progress in Mathematics, 266. Birkh\"auser Verlag, Basel, 2008.
\end{thebibliography}
\end{document}